\documentclass[12pt]{amsart}
\usepackage[utf8]{inputenc}
\usepackage{amsmath,amsthm}
\usepackage{multirow}
\usepackage{graphics}
\usepackage{tikz}
\usetikzlibrary{arrows.meta}
\usepackage{listings}

\newtheorem{theorem}{Theorem}

\begin{document}

\title[Scriabin, the mystic chord and its counterpoint]{Contrapuntal Aspects of the Mystic Chord and Scriabin's Piano Sonata No. 5}
\author{Octavio A. Agustín-Aquino \and Guerino Mazzola}
\address{Instituto de Física y Matemáticas, Universidad Tecnológica de la Mixteca, Huajuapan de León, Oaxaca, México}
\email{octavioalberto@mixteco.utm.mx}
\address{School of Music, University of Minnesota, MN, USA}
\email{mazzola@umn.edu}

\date{September 25, 2018}
\abstract{
We present statistical evidence for the importance of
the ``mystic chord'' in Scriabin's Piano Sonata No. 5, Op. 53, from
a computational and mathematical counterpoint perspective. More
specifically, we compute the effect sizes and $\chi^{2}$ tests with respect to
the distributions of counterpoint symmetries in the Fuxian, mystic,
Ionian and representatives of the other three possible counterpoint
worlds in two passages of the work, which provide evidence of a
qualitative change between the Fuxian and the mystic world in the sonata.}
\thanks{This work was partially supported by a grant from the \emph{Niels Hendrik Abel Board}.}}
\maketitle

\section{Introduction}
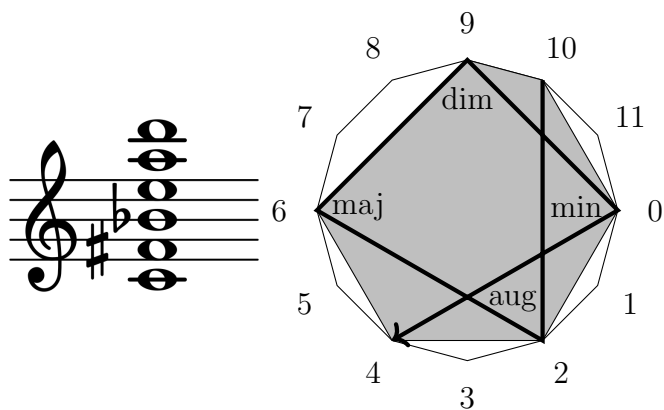
\begin{figure}[ht]
\begin{center}
\raisebox{4em}{\includegraphics[scale=0.2]{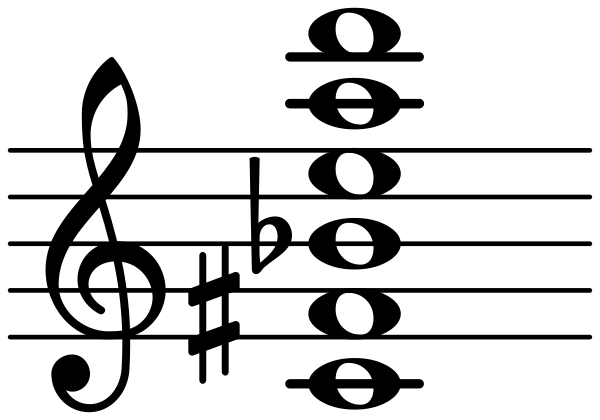}}%
\begin{tikzpicture}
     \draw    (2.000000000000000,-0.000000000000000)
--(1.732050807568877,-1.000000000000000)
--(1.000000000000000,-1.732050807568877)
--(0.000000000000000,-2.000000000000000)
--(-1.000000000000000,-1.732050807568877)
--(-1.732050807568877,-1.000000000000000)
--(-2.000000000000000,-0.000000000000000)
--(-1.732050807568878, 0.999999999999999)
--(-1.000000000000001, 1.732050807568877)
--(-0.000000000000000, 2.000000000000000)
--(1.000000000000000, 1.732050807568877)
--(1.732050807568877, 1.000000000000001)
--(2.000000000000000,-0.000000000000000);
   \node at (2.500000000000000, -0.000000000000000) {$0$};
   \node at (2.165063509461097, -1.250000000000000) {$1$};
   \node at (1.250000000000000, -2.165063509461096) {$2$};
   \node at (0.000000000000000, -2.500000000000000) {$3$};
  \node at (-1.250000000000000, -2.165063509461097) {$4$};
  \node at (-2.165063509461097, -1.250000000000000) {$5$};
  \node at (-2.500000000000000, -0.000000000000000) {$6$};
  \node at (-2.165063509461097, 1.249999999999999) {$7$};
  \node at (-1.250000000000001, 2.165063509461096) {$8$};
  \node at (-0.000000000000000, 2.500000000000000) {$9$};
   \node at (1.250000000000000, 2.165063509461096) {$10$};
   \node at (2.165063509461096, 1.250000000000001) {$11$};
   
   \draw [fill=lightgray](2.000000000000000,-0.000000000000000)--(1.000000000000000,-1.732050807568877)--(-1.000000000000000,-1.732050807568877)--(-2.000000000000000,-0.000000000000000)--(-0.000000000000000, 2.000000000000000)--(1.000000000000000, 1.732050807568877)--(2.000000000000000,-0.000000000000000);
   \draw [ultra thick,->] (1.000000000000000, 1.732050807568877)--(1.000000000000000,-1.732050807568877)--(-2.000000000000000,-0.000000000000000)--(-0.000000000000000, 2.000000000000000)--(2.000000000000000,-0.000000000000000)--(-1.000000000000000,-1.732050807568877);
   
   \node at (0.600000000000000,-1.209038105676658) {aug};
   \node at (-1.450000000000000, 0.050000000000000) {maj};
   \node at (-0.000000000000000, 1.500000000000000) {dim};
   \node at (1.450000000000000, 0.050000000000000) {min};
    \end{tikzpicture}
\end{center}
\caption{Scriabin's mystic chord and four prominent triads (augmented, major, diminished, and minor) which cover it.}
\label{F:Mystic}
\end{figure}

A prominent chord in Alexander Scriabin's late work is the so-called ``mystic chord'' or ``Prometheus chord'', whose pc-set
when the root is C is $M=\{0,6,10,4,9,2\}$ \cite[p. 23]{vD14}.  It can be seen as a chain of thirds, and thus can be covered by an
augmented triad followed by a diminished triad, together with a major triad followed by a minor triad (see Figure
\ref{F:Mystic}).
This surely evidences the strong tonal ambiguity of the chord, which is also associated to the Impressionism in music during
the late 19th and early 20th centuries in Europe \cite{rT99}. The mystic chord can be seen as an extension of the French sixth, which is completely
contained in a whole tone scale; the mystic chord also has this property safe for a ``sensible'' tone \cite[p. 278]{BBB17}. As we will see, this is
an important feature from the perspective of the mathematical counterpoint theory developed by Mazzola \cite[Part VII]{gM17}.

With respect to the structural role of the mystic chord in Scriabin's works, Gottfried Eberle states\footnote{All the musical events of the work, harmonic
as well as melodic and contrapuntal, are essentially within this six-tone complex: ``[...] all the melodic voices are on the sounds of the accompanying harmony, all counterpoints are subordinated to the same principle''.}, for instance \cite[p. 14]{gE78},

\begin{quote}
 Alle musikalischen Ereignisse des Werks, harmonische wie melodische und kontrapunktische, seien im wesentlichen
in diesem sechstönigen Komplex gegründet: ``[...] alle melodischen Stimmen sind auf den Klängen der
begleitenden Harmonie gebaut, alle Kontrapunkte sind demselben Prinzip untergeordnet'' [Sabanajew, 1912].
\end{quote}

and, moreover\footnote{The statement is correct: The basic harmony of Prometheus is no longer understood and treated by Scriabin as dissonance:
``This is a basic harmony, a consonance''.} \cite[p. 16]{gE78},

\begin{quote}
Die Aussage trifft zu: Die Grundharmonie des Prometheus wird von Skrj\-abin nicht läger als Dissonanz begriffen und behandelt: ``Das ist eine Grundharmonie, eine Konsonanz'' [Sabaneev, 1925].
\end{quote}

More specifically, in our contrapuntal interpretation, we view the mystic chord defining what is called a \emph{strong dichotomy}, i. e., a bipartition of
the pitch class set $\mathbb{Z}_{12}$
such that its only affine symmetry is the identity \cite[Part VII, Chapter 30]{gM17}. In particular, it\footnote{It is
interesting to note that the mystic chord can be seen as the first of a sequence of
dichotomies that are always strong in microtonal equitempered tunings, as described in \cite[Proposition 2.4 and Remark 3]{oA09}.}
belongs to the class 78 in Mazzola's classification as it
appears in his \emph{Topos of Music} \cite[Table L.1]{gM17}. A strong dichotomy can be understood as a division of pitch classes in generalized consonances and
dissonances, because the classical consonances of Renaissance counterpoint also define a strong dichotomy. There are four other bipartitions with the aforementioned mathematical properties (modulo affine symmetries). Here
is a list of selected representatives of the strong dichotomies:
\begin{align*}
(K/D) &= \Delta_{82} = (\{0, 3, 4, 7, 8, 9\}/\{1, 2, 5, 6, 10, 11\}),\\
(M/N ) &= \Delta_{78} = (\{0, 2, 4, 6, 9, 10\}/\{1, 3, 5, 7, 8, 11\}),\\
(I/J) &= \Delta_{64} = (\{2, 4, 5, 7, 9, 11\}/\{0, 1, 3, 6, 8, 10\}),\\
\Delta_{68} &= (\{0, 1, 2, 3, 5, 8\}/\{4, 6, 7, 9, 10, 11\}),\\
\Delta_{71} &= (\{0, 1, 2, 3, 6, 7\}/\{4, 5, 8, 9, 10, 11\}),\\
\Delta_{75} &= (\{0, 1, 2, 4, 6, 8\}/\{3, 6, 7, 9, 10, 11\}).
\end{align*}

The subindex represents the class number in Mazzola's classification, whereas $(K/D)$, $(M/N)$ and $(I/J)$ are the \emph{Fuxian}, \emph{mystic} and
\emph{Ionian}\footnote{This name stems from the fact that this representative consists in all proper (non-vanishing) intervals in the Ionian mode,
when counted from the tonic.} dichotomies, respectively.

If we study the predicted allowed steps for a counterpoint distilled from $M$, we find that a particularly favorable scale for cantus firmus
pitches is one particular transposed mode of the whole-tone scale, namely the one with pc-set $\{1,3,5,7,9,11\}$, with only eight forbidden transitions
if we do not mind if the discantus leaves the selected whole-tone scale, or four if the discantus has to remain within the scale. Thus, we may consider two
representatives of mystic chord: one, like $M$, which shares most
of its tones with the ``even'' whole-tone scale $\{0,2,4,6,8,10\}$, and the other\footnote{We denote with $T^{x}$
the transposition by $x$.} one $T^{1}M$ that is closer to the ``odd'' whole-tone scale.

Hence, in general, for the even mystic chord, a very good scale for counterpoint is the odd whole-tone scale, and vice versa.

Unfortunately, for the whole-tone scale there is no analogue of Noll's theorem connecting a harmony based on triads and counterpoint (as
explained in \cite[Section 30.2.1]{gM17}), since among all possible triads there is none whose set of endomorphisms is such that their linear part yields a strong dichotomy. This, by the way, is in accordance
with classical musicological opinion on the scale of its poor harmonic possibilities (at least from the tonal harmony perspective \cite[p. 486]{KP95}), and perhaps it was an attractive characteristic for Scriabin to use it in his music.

\section{A Quick Overview of Mazzola's Counterpoint Model}

Before we proceed, let us make a remark on
notation: we denote a \emph{counterpoint interval} by $(x,y)$, where $x$ is the \emph{cantus firmus} and $y$ is the interval that separates
it from the \emph{discantus}. Thus, $(2,7)$ represents a counterpoint interval where the cantus firmus is $2$ and
the discantus is $9$, because the separation between them, modulo octaves, is $(9-2)\bmod 12 = 7$.

In Mazzola's counterpoint model all the pitches are considered modulo octave. Thus the intervals between two tones
reduces to $\mathbb{Z}_{12}$. In particular, as far as Renaissance counterpoint and the famous
Fux's treatise \emph{Gradus ad Parnassum} \cite{jF65} are concerned, the set of consonances is $K=\{0,3,4,7,8,9\}$ and thus dissonances are $D=\mathbb{Z}_{12}\setminus K$. The bipartition of intervals $(K/D)$ is an example of a \emph{strong dichotomy}, which
we shall define now. First, the group of affine symmetries between pitch classes in $\mathbb{Z}_{12}$ consists of those of the form
\[
 T^{u}.v(x) = vx+u
\]
where $v$ is coprime with $12$, i. e., $v=1,5,7,11$ and $u\in\mathbb{Z}_{12}$. Note that the affine symmetry
$p=T^{2}.5$ is such that $p(K)=D$ (acting pointwise) and it is the only one with this property. It is called
the \emph{polarity} of the set of consonances. Precisely those dichotomies that possess a unique polarity are called \emph{strong}. 
As we have already mentioned, there is a total of six strong dichotomies with these properties up to equivalence
under the action of the group of affine symmetries.

Counterpoint intervals can be endowed with the structure of
dual numbers\footnote{The pairs $(x,y)$ can also be written as $x+\epsilon.y$ with $\epsilon^{2}=0$ in
commutative algebra. The reason to introduce this algebraic structure is that it describes tangent vectors (see \cite[Section 7.5]{gM17} for further details).} $\mathbb{Z}_{12}[\epsilon]\in\mathbb{Z}_{12}\times \mathbb{Z}_{12}$, defining the sum
\[
(x_{1},y_{1})+(x_{2},y_{2}) = (x_{1}+x_{2},y_{1}+y_{2})
\]
and the multiplication
\begin{equation}\label{E:DefMult}
(x_{1},y_{1})(x_{2},y_{2}) = (x_{1}x_{2},x_{1}y_{2}+x_{2}y_{1}).
\end{equation}

The group of symmetries for counterpoint intervals consists of symmetries of the form
\[
 T^{(u_{1},u_{2})}.(v_{1},v_{2})
\]
with $(v_{1},v_{2})$ an invertible element with respect to the multiplication defined by \eqref{E:DefMult}, which
amounts to require $v_{1}$ to be invertible. We denote them with $\overrightarrow{GL}(\mathbb{Z}_{12}[\epsilon])$. For an arbitrary strong dichotomy
$(X/Y)$ such that its set of consonances and dissonances are $X$ and $Y$, respectively, the set of all consonant intervals is
\[
X[\epsilon]:=\{(c,x):c\in\mathbb{Z}_{12},x\in X\}
\]
and the set of dissonant intervals is $Y[\epsilon] := \mathbb{Z}_{12}\times \mathbb{Z}_{12}\setminus X[\epsilon]$.
In particular, there exists a canonical symmetry $p_{c}$
such that $p_{c}(X[\epsilon]) = Y[\epsilon]$ and leaves the intervals with cantus firmus $c$ invariant, and it is called an \emph{induced polarity}.
We let the group
$\overrightarrow{GL}(\mathbb{Z}_{12}[\epsilon])$ to act pointwise on subsets $S\subseteq \mathbb{Z}_{12}[\epsilon]$, and we call
\[
gX[\epsilon] = \{g(c,x):(c,x)\in X[\epsilon]\}
\]
a set of \emph{$g$-deformed consonant intervals}. The \emph{$g$-deformed dissonant intervals} are, of course, $gY[\epsilon]$.

Now a \emph{counterpoint symmetry} $g$ for a consonant interval $\xi = (c,x)$ is one such that
\begin{enumerate}
\item the interval $\xi$ is a $g$-deformed dissonant interval,
\item it commutes with the induced polarity $p_{x}$, that is, $p_{x}\circ g = g\circ p_{x}$, and
\item the set of consonances that are also $g$-deformed consonant intervals is as large as possible within the symmetries with the above properties.
\end{enumerate}

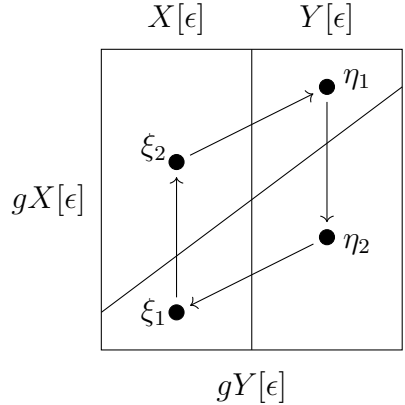
\begin{figure}[ht]
\begin{center}
    \begin{tikzpicture}
        \draw (0,0) rectangle (4,4);
		\draw (2,0)--(2,4);
		\draw (0,0.5) -- (4,3.5);
        \draw [fill] (1,0.5) circle [radius=0.1];
		\draw [fill] (1,2.5) circle [radius=0.1];
		\draw [fill] (3,1.5) circle [radius=0.1];
		\draw [fill] (3,3.5) circle [radius=0.1];
		\draw [->] (1,0.7)--(1,2.3);
		\draw [->] (3,3.3)--(3,1.7);
		\draw [->] (2.82,1.41)--(1.18,0.59);
		\draw [->] (1.18,2.59)--(2.82,3.41);
		\node at (1,4.4) {$X[\epsilon]$};
		\node at (3.0,4.4) {$Y[\epsilon]$};
		\node at (-0.7,2) {$gX[\epsilon]$};
		\node at (2,-0.5) {$gY[\epsilon]$};
		\node at (0.7,0.5) {$\xi_{1}$};
		\node at (0.7,2.7) {$\xi_{2}$};
		\node at (3.4,3.6) {$\eta_{1}$};
		\node at (3.4,1.4) {$\eta_{2}$};
    \end{tikzpicture}
\end{center}	
\caption{Deformed consonant and dissonant intervals, and possible transitions.}
\label{F:Deformation}
\end{figure}

The consonant intervals which are $g$-deformed consonant intervals simultaneously for a counterpoint symmetry $g$ are called the
\emph{admitted successors}. The idea behind these definitions is that transitions from consonance to consonance do not exhibit tension explicitly. Mazzola's solution to reveal this concealed tension is to see the first consonance as
a deformed dissonance, and then resolving it to a deformed consonance that it is also a consonance
(the step from $\xi_{1}$ to $\xi_{2}$ in Figure \ref{F:Deformation}).

It is important to mention that counterpoint symmetries can be calculated for cantus firmus $c=0$ and then
translated suitably \cite[Section 31.3.1]{gM17}. It is also relatively straightforward to adapt the calculations' results for other consonances
and dissonances that are affinely equivalent to the $(X/Y)$ dichotomy without redoing them entirely \cite[Section 31.3.4]{gM17}.
We should stress that, when $(X/Y)$ is the Fuxian dichotomy, then we recover many salient features of Renaissance counterpoint
theory via counterpoint symmetries and admitted successors, for example: the fourth is classified as a dissonance, there is a general prohibition
of parallel fifths, the tritone rules hold in the so-called reduced strict style (which is obtained to applying Fux's rules modulo octave), and
the major scale is optimal for contrapuntal purposes, in the sense that it allows the largest number of allowed steps. See \cite[Chapter 3]{AJM15}
and \cite[Section 31.4.1]{gM17} for further details.

\section{General Transitions in Mazzola's Counterpoint Model}

The aforementioned model apparently only handles the transitions from consonances to consonances. But upon
reflection we realize that it also handles the dissonance-to-dissonance steps (like the one from
$\eta_{1}$ to $\eta_{2}$ in Figure \ref{F:Deformation}), by simply applying the
induced polarity $p_{0}$ and translating the results accordingly\footnote{This, by the way, leads to a
natural concept of dissonant counterpoint \cite{hC96}.}. More explicitly, the admitted successors of a dissonant interval $\eta_{1}\in D[\epsilon]$ for a given counterpoint symmetry $g$ are translates of
\begin{align*}
p_{0}(gK[\epsilon]\cap K[\epsilon])&= p_{0}gK[\epsilon]\cap p_{0} K[\epsilon]\\
&= gp_{0}K[\epsilon]\cap D[\epsilon] = gD[\epsilon]\cap D[\epsilon].
\end{align*}

Nevertheless, it is less obvious that
dissonance-to-consonance steps or \emph{resolutions} \cite[p. 56]{jF65} are also modeled\footnote{In Renaissance counterpoint the notion of
resolution is understood only by stepwise movement of voices \cite[p. 131]{kJ92}, but the model can trivially be restricted to fulfill this
requirement.}, like the one from $\eta_{2}$ to $\xi_{1}$ in Figure \ref{F:Deformation}. In fact, we ought to define a ``crossing''
counterpoint symmetry $g$ for a dissonance $\eta_{2}$ as one such that, apart from the obvious requirement of
commuting with the induced polarity, $\eta_{2}$ is an appropriate deformation of a dissonant interval and maximizes the intersection
of the original consonant intervals and $g$-deformed dissonant intervals. But the following result shows that no true generalization is needed.

\begin{theorem} Let $\eta\in Y[\epsilon]$. The symmetry $g'\in \overrightarrow{GL}(\mathbb{Z}_{12}[\epsilon])$
satisfies that
\begin{enumerate}
\item the interval $\eta$ belongs to $g'Y[\epsilon]$,
\item it commutes with $p_{0}$ and
\item the set $g'Y[\epsilon]\cap X[\epsilon]$ has the largest cardinality among the symmetries with the
previous two properties
\end{enumerate}
if, and only if, $g=p_{0}\circ g'$ is a counterpoint symmetry for the interval $\xi=p_{0}(\eta)\in X[\epsilon]$.
\end{theorem}
\begin{proof}
Note first that $p_{0}$ is involutive, i. e., $p_{0}=p_{0}^{-1}$. Thus $g'$ commutes with $p_{0}$ if and
only if $g=p_{0}\circ g'$ does, since
\[
 g=p_{0}\circ g' = g'\circ p_{0}\iff g'=p_{0}\circ g'\circ p_{0} = g\circ p_{0} = p_{0}\circ g.
\]

Thus
\[
\xi \in gY[\epsilon] \iff \eta = p_{0}(\xi)\in p_{0}\circ g Y[\epsilon] = p_{0}\circ p_{0}\circ g' Y[\epsilon] = g'Y[\epsilon].
\]

Finally, because
\begin{align*}
 gX[\epsilon]\cap X[\epsilon] &= g\circ p_{0} Y[\epsilon]\cap X[\epsilon]\\
 &=g'Y[\epsilon]\cap X[\epsilon],
\end{align*}
the maximization of the cardinality of the admitted successors $gX[\epsilon]\cap X[\epsilon]$ is equivalent to
the maximization of $g'Y[\epsilon]\cap X[\epsilon]$.
\end{proof}

In other words: in order to model the transition from dissonance to consonance with symmetries, we regard the dissonance $\eta$ as a $g$-deformed
dissonant interval and we admit as a successor a deformed dissonant interval that it is also a consonant interval, hence \emph{dissolving} the contrast between
dissonance and consonance via a deformation. The case for consonance-to-dissonance steps or \emph{preparations} \cite[p. 44]{tB05}, like the
one from $\xi_{2}$ to $\eta_{1}$ in Figure \ref{F:Deformation}, is not only analogous: it is symmetrical, therefore we omit the details here.

\section{Some Statistical Contrapuntal Properties of the Counterpoint Worlds}\label{Sec:MundosEstadistica}

In Table \ref{T:DistFuxMystic} we see the distribution of the number of contrapuntal symmetries between all intervals for the selected representatives
mentioned in Section 1 of the \emph{counterpoint worlds}\footnote{A counterpoint world is a directed graph, where each vertex is a counterpoint interval and there is an arrow connecting for each valid step. See \cite[Chapter 4]{AJM15} for further details.}.

If you have, for instance, the step $((0,3),(2,4))$ in the Fuxian world, there
are two counterpoint symmetries that ``allow'' it. Contrariwise, the parallel fifths step $((0,7),(2,7))$ is forbidden, for it has
$0$ counterpoint symmetries. In the Fuxian world the maximum number of symmetries mediating in a step is $5$.
For the $12^{4}$ possible steps
(besides consonance to consonance also dissonance to dissonance, dissonance to consonance, and consonance to dissonance steps can be considered in the
 extended model),
we have a total of $2400$ inadmissible steps, $4992$ steps with only one counterpoint symmetry, and so on within the Fuxian world.
In Table \ref{T:StatsFuxMystic} we find the mean and standard deviation of the distributions of Table \ref{T:DistFuxMystic}, as well as the probability of
a random pair of intervals to be valid.

It should be noted that the Fuxian world has the greatest liberty for counterpoint transitions, since it has the minimum of forbidden steps.
It is followed by the $\Delta_{68}$, mystic, Ionian, $\Delta_{75}$, and $\Delta_{71}$, in that order.

Let $p_{\Delta_{1}}$, $p_{\Delta_{2}}$, and $p_{\Delta_{1}\cap\Delta_{2}}$ the
probability of a random step to valid in worlds $\Delta_{1}$, $\Delta_{2}$, and
in both, respectively. The absolute value of
$|p_{\Delta_{1}} p_{\Delta_{2}}-p_{\Delta_{1}\cap\Delta_{2}} |$, which measures the deviation of the events from independence, appears in Table \ref{T:Independencia}. Notice that the highest deviation from independence occurs between the Ionian and
$\Delta_{75}$ worlds. Nevertheless, most of the other pairs deviate from
independence approximately by a 5\% of this maximum value or even less.

\begin{table}[ht]
\begin{center}\begin{tabular}{|c|c|c|c|c|c|c|}
\hline
Number of cpt. & Fuxian  & Mystic & Ionian & $\Delta_{68}$ & $\Delta_{71}$ & $\Delta_{75}$\\
symmetries & steps &steps & steps & steps & steps &steps\\
\hline
$0$&$2400$&$3840$	&	$4224$	&	$3744$	&	$4608$	&	$4320$\\
$1$&$4992$&$7296$	&	$11136$	&	$6912$	&	$6912$	&	$9120$\\
$2$&$9120$&$6720$	&	$4608$	&	$6624$	&	$6336$	&	$4992$\\
$3$&$384$&$0$		&	$768$	&	$3456$	&	$0$		&	$1440$\\
$4$&$2304$&$2880$	&	$0$		&	$0$		&	$2880$	&	$864$\\
$5$&$1536$&$0$		&	$0$		&	$0$		&	$0$		&	$0$\\
\hline
\end{tabular}\end{center}

\caption{Distribution of the number of counterpoint symmetries for selected representatives of the counterpoint worlds.}
\label{T:DistFuxMystic}
\end{table}

\begin{table}[ht]
\begin{center}\begin{tabular}{|c|c|c|c|}
\hline
Counterpoint world & Mean & Std. dev. & Prob. of admissibility\\
\hline
Fuxian & $1.99074$ & $1.35401$ & $0.88426$\\
Mystic  & $1.55556$ & $1.20444$ & $0.81481$\\
Ionian & $1.09259$ & $0.75202$ & $0.79630$\\
$\Delta_{68}$ & $1.47222$ & $0.97145$ & $0.81944$\\
$\Delta_{71}$ & $1.50000$ & $1.23606$ & $0.77778$\\
$\Delta_{75}$ & $1.29630$ & $1.00703$ & $0.79167$\\
\hline
\end{tabular}\end{center}

\caption{Mean and standard deviation of the number
of counterpoint symmetries in steps within the counterpoint worlds, and the
probability that two random intervals to be a valid succession.}
\label{T:StatsFuxMystic}
\end{table}

\begin{table}[ht]
\begin{center}\begin{tabular}{|c|c|c|c|c|c|c|}
\hline
       &Fuxian &Mystic &Ionian &$\Delta_{68}$ &$\Delta_{71}$ &$\Delta_{75}$\\
\hline
Fuxian 			&---	&0.0303	&0.0154	&0.0231	&0.0113	&0.0129\\
\hline
Mystic 			&		&---	&0.0564	&0.0036	&0.0190	&0.0000\\
\hline
Ionian 			&		&		&---	&0.0067	&0.6277	&0.6373\\
\hline
$\Delta_{68}$   &		&		&		&—--	&0.0050	&0.0006\\
\hline
$\Delta_{71}$	&		&		&		&		&---	&0.0046\\
\hline
\end{tabular}\end{center}
\caption{Absolute value of the difference between the
product of the probabilities of being valid for a random
step in two worlds and the probability of being valid in
both worlds.}
\label{T:Independencia}
\end{table}

\section{Fuxian and Mystic worlds in Scriabin's Piano Sonata No. 5, op. 53}

We can find some explanatory power of Mazzola's contrapuntal model with Scriabin's Piano Sonata No. 5, op. 53 \cite{aS07}, which is notable for the
explicitness of the mystic chord.

In the prologue section \cite[Chapter IV]{hW87}, which we will call part 1, spanning measures 13 to 46, taking in most of the cases the cantus firmus as E
(as suggested by standard musicological analysis of the work \cite[pp. 2-3]{sZ10}), we have $37$ contrapuntal transitions within measures 13--31.
Next we isolate 36 possible transitions (omitting repetitions of certain patterns) for measures from 47 to 61, which are the initial measures of what is known as the first exposition of the sonata \cite[Chapter IV]{hW87},
and that we will call part 2. The average number of counterpoint
symmetries per step and the standard deviation appears in Table \ref{T:AvgNStdDev}.

\begin{table}[ht]
\begin{center}\begin{tabular}{|c|c|c|c|c|}
\hline
\multirow{2}{*}{Cpt. world} & \multicolumn{2}{l|}{Part 1} & \multicolumn{2}{l|}{Part 2} \\ \cline{2-5} 
                            & Avg. \# of sym. & Std. dev. & Avg. \# of sym. & Std. dev. \\ \hline
Fuxian			&$2.10811$	&$1.10010$	&$1.94444$	&$0.62994$\\
\hline
Mystic			&$2.16216$	&$1.34399$	&$1.38889$	&$1.20185$\\
\hline
Ionian			&$1.02703$	&$0.89711$	&$0.88889$	&$0.74748$\\
\hline
$\Delta_{68}$	&$1.54054$	&$0.90045$	&$1.77778$	&$0.89797$\\
\hline
$\Delta_{71}$	&$1.59459$	&$1.32202$	&$1.27778$	&$0.97427$\\
\hline
$\Delta_{75}$	&$1.13514$	&$0.75138$	&$1.08333$	&$0.69179$\\
\hline
\end{tabular}\end{center}
\caption{Average number of symmetries per step and
standard deviation in part 1 (measures 13–31) and part
2 (measures 47–61) of Scriabin's Sonata.}
\label{T:AvgNStdDev}
\end{table}

\begin{table}[ht]
\begin{center}\begin{tabular}{|c|c|c|c|c|}
\hline
\multirow{2}{*}{Cpt. world} & \multicolumn{2}{l|}{Part 1} & \multicolumn{2}{l|}{Part 2} \\ \cline{2-5} 
                            & ES & 95\% conf. intvl. & ES & 95\% conf. intvl. \\ \hline
Fuxian			&$0.087$	&$[-0.236, 0.409]$	&$-0.034$	&$[-0.361, 0.293]$\\
\hline
Mystic			&$0.504$	&$[0.181, 0.826]$	&$-0.138$	&$[-0.465, 0.189]$\\
\hline
Ionian			&$-0.087$	&$[-0.410, 0.235]$	&$-0.271$	&$[-0.598, -0.056]$\\
\hline
$\Delta_{68}$	&$0.070$	&$[-0.252,0.393]$	&$0.315$	&$[-0.012, 0.642]$\\
\hline
$\Delta_{71}$	&$0.077$	&$[-0.246, 0.399]$	&$-0.180$	&$[-0.507, 0.147]$\\
\hline
$\Delta_{75}$	&$-0.160$	&$[-0.483, 0.162]$	&$-0.212$	&$[-0.539, 0.115]$\\
\hline
\end{tabular}\end{center}
\caption{Effect size and 95\% confidence intervals for
part 1 (measures 13–31) and part 2 (measures 47–61) of
Scriabin's Sonata.}
\label{T:EffSize}
\end{table}

We now compute the effect sizes\footnote{The effect size we take is the so-called Cohen's $d$, which is the mean
difference on the means between the two variables divided by the pooled standard deviation. See \cite{gC12} for further details.} (ES) and the corresponding 95\% confidence intervals for the two parts of the work under analysis (see Table
\ref{T:EffSize}), and we clearly observe a relatively large positive effect of the mystic world in the first
part. For the second part, the presence of all the worlds reduces, except for the
$\Delta_{68}$ world, which is the only one whose effect size \emph{increases}. This strange phenomenon,
along with the apparent absence of the Fuxian world, is further clarified by the following tests.

\begin{table}[ht]
\begin{center}\begin{tabular}{|l|c|c|}
\hline
\multirow{2}{*}{Counterpoint world} & \multicolumn{2}{l|}{$p$-value} \\ \cline{2-3} 
                                    & Part 1         & Part 2        \\ \hline
Fuxian &$3.1778\times 10^{-9}$&$ 2.2141\times 10^{-3}$\\
\hline
Mystic &$0.0219$&$0.0436$\\
\hline
Ionian &$0.2084$&$0.1809$\\
\hline
$\Delta_{68}$ &$0.6460$&$0.2669$\\
\hline
$\Delta_{71}$ &$0.7120$&$0.3535$\\
\hline
$\Delta_{75}$ &$0.5320$&$0.3007$\\
\hline
\end{tabular}\end{center}
\caption{List of $p$-values for $\chi^{2}$ tests for part 1 (measures 13–31) and part 2 (measures 47–61) of Scriabin's
Sonata.}
\label{T:Chitests}
\end{table}

If we perform $\chi^{2}$ tests \cite[Chapter 8]{rW12} for the frequencies of number of symmetries, with
$p$-values appearing in Table \ref{T:Chitests}, we see that the only worlds that would pass a 95\% test
would be the distributions of the Fuxian and mystic worlds, which confirm that their presences are the only
significant ones. If we now restrict the $\chi^{2}$ test to the permitted versus allowed steps frequencies\footnote{In
fact, the kurtosis of the distribution of the number of symmetries per step in the Fuxian world are $4.70239$ and $7.55462$ for
part 1 and 2, respectively. This means that the second distribution deviates less from its mean, and thus in this case
Cohen's $d$ does not explain the change sufficiently because the mean of both distributions is very close to
the general one. This was also observed in the
first-species fragments of \emph{Misae Papae Marcelli} by G. P.
Palestrina against the general distribution; see \cite{AN10} for details.}, we find
the values in Table \ref{T:Chitests2}. Now we notice that the mystic world would pass a 89\% test for part 1 and not
for part 2, and the converse is true for the Fuxian world.

\begin{table}[ht]
\begin{center}
\begin{tabular}{|l|c|c|c|c|}
\hline
\multirow{3}{*}{Cpt. world} & \multicolumn{2}{l|}{Part 1}     & \multicolumn{2}{l|}{Part 2} \\ \cline{2-5}
  &Fraction of     & $\chi^{2}$ & Fraction of     & $\chi^{2}$   \\
        & permitted steps & statistic  & permitted steps & statistic\\
\hline
Fuxian             &$0.92$ & $0.5099$ & $1.00$ & $0.0300$ \\ \hline
Mystic             &$0.91$ & $0.1031$ & $0.83$ & $0.7748$ \\ \hline
\end{tabular}
\end{center}
\caption{List of fractions of permitted steps per part and $p$-values for $\chi^{2}$ tests for frequencies of permitted/forbidden steps in part 1 (measures 13–31) and part 2 (measures 47–61) of Scriabin's Sonata for the Fuxian and mystic worlds.}
\label{T:Chitests2}
\end{table}

In other words: while
part 2 does not use functional tonal harmony, it is much closer to the counterpoint
of the standard consonances than to one stemming from the mystic chord heard in part 1, and no
other counterpoint worlds are evident aside from these.

\begin{table}[ht]
\begin{center}\begin{tabular}{|c|c|c|}
\hline
Counterpoint 	&Consonances in part 1	&Consonances in part 2\\
world			&(57)					&(50)\\
\hline
Fuxian			&$24$					&$44$\\
\hline
Mystic			&$39$					&$19$\\
\hline
Ionian			&$31$					&$32$\\
\hline
$\Delta_{68}$	&$17$					&$24$\\
\hline
$\Delta_{71}$	&$22$					&$17$\\
\hline
$\Delta_{75}$	&$24$					&$23$\\
\hline
\end{tabular}\end{center}
\caption{Total number of consonances in each of the analyzed parts of Scriabin's sonata (total number
of intervals is in parentheses).}
\label{T:TotalesConsonancias}
\end{table}

\begin{table}
\begin{center}\begin{tabular}{|c|l|l|l|l|l|l|l|l|}
\hline
\multirow{2}{*}{Counterpoint world} & \multicolumn{4}{l|}{Part 1} & \multicolumn{4}{l|}{Part 2} \\ \cline{2-9} 
                                    & DD    & DC    & CD   & CC   & DD    & DC    & CD   & CC   \\ \hline
Fuxian					&$14$		&$8$		&$6$		&$9$		&$4$		&$0$		&$0$		&$32$\\
\hline
Mystic					&$1$		&$8$		&$13$		&$15$		&$6$		&$18$		&$12$		&$0$\\
\hline
Ionian					&$5$		&$12$		&$11$		&$9$		&$0$		&$15$		&$9$		&$12$\\
\hline
$\Delta_{68}$			&$17$		&$11$		&$7$		&$2$		&$7$		&$10$		&$16$		&$3$\\
\hline
$\Delta_{71}$			&$15$		&$8$		&$8$		&$6$		&$13$		&$10$		&$13$		&$0$\\
\hline
$\Delta_{75}$			&$22$		&$0$		&$0$		&$15$		&$19$		&$0$		&$0$		&$17$\\
\hline
\end{tabular}\end{center}
\caption{Total number of different transitions (DD: dissonance to dissonance; DC: dissonance to consonance; CD: consonance to dissonance;
DD: dissonance to dissonance) for each of the analyzed parts of Scriabin's sonata.}
\label{T:TotalesTransiciones}
\end{table}

\section{Some Additional Observations}

Although it is not directly connected to Mazzola's counterpoint theory (since it only depends on the
purely combinatorial classification of consonances and dissonances), it is very interesting to note
that the highest count of consonances for part 1 under analysis corresponds to the mystic world, whereas
the maximum occurs for the Fuxian consonances in part 2. If we calculate the number of transitions
for the four possible transitions between consonances and dissonances, we note that the maximum number
of consonance-to-consonance steps in part 1 occurs for the mystic and $\Delta_{75}$ worlds (although
the display of resolutions and preparations for the mystic world is evident and totally absent for the
world $\Delta_{75}$). The maximum number of consonant transitions goes for the Fuxian world in part 2,
as expected.

While not as explicit as an extraction of first species counterpoint from the sonata, we can
find more evidences of the importance of the mystic chord as a choice of consonances and the
role the whole tone scale has with respect to it.

For instance, from measure 102 to 103 Scriabin favors pitches within the even whole-tone scale and in 104 and 105 he states an odd mystic chord;
then he suddenly changes the key and begins to stress the even whole-tone scale.

Another similar situation occurs in measures 130 and 131, where Scriabin displays an arpeggiated even mystic chord followed by an arpeggiated
odd whole-tone scale in the following measure. Quite interestingly, this is continued by an odd mystic chord in measures 136 and 137, but preceding
it with and ambiguity between the even and odd whole tone scale, anticipating another sudden change of mood in measure 140.

A final explicit apparition of an even mystic chord in measure 262 is also associated with a dynamical fluctuation in the piece, but in this case its interaction
with the whole scale is less apparent but seems to be in favor of the odd whole-tone scale, as expected.

\section{Conclusions}

As Eberle and other scholars who specialize in Alexander Scriabin have pointed out, the mystic or Prometheus chord has been a key architectural principle in his works but its relation to the contrapuntal aspect of them has largely been neglected or not understood.
Through an extension of Mazzola's counterpoint model, where the mystic chord can literally (as Scriabin himself claimed) be taken as the consonances, a counterpoint theory emerges such that general transitions between consonances and dissonances can be handled, and thus we can compare the contrapuntal content of two different passages of Scriabin's fifth piano sonata not only across one but all of the counterpoint worlds. The fact that we can perceive Scriabin's accomplishment of the combination of two counterpoint worlds within one composition attests the power of the mathematical model for understanding difficult works of art and project them into the future.

\section*{Acknowledgments}

We thank Thomas Noll at Escola Superior de Música de Catalunya and Daniel Tompkins at Florida State University for their
valuable feedback.

\section*{Appendix: Source Code and Data}

The following code implements a function in Octave (version 4.2.0) to calculate the number of counterpoint
symmetries per step in a sequence of counterpoint intervals, encoded as columns of a matrix.

\begin{lstlisting}[language={Octave},basicstyle={\footnotesize}]
function R = analisis(M, K, simetrias)

% --- R = analisis (M, K, simetrias)
% Calculates the number of counterpoint symmetries
% that mediate between the counterpoint steps in M,
% where the first row corresponds to the cantus
% firmus and the second row to the intervals with
% the discantus. If the parameters K and simetrias
% are provided, it uses K as the consonances, and
% within the array simetrias they should appear in
% the format [a b c], corresponding to
% e^(0,c).(a,b).

if(nargin()==1)

 K = [0 3 4 7 8 9];

 simetrias = {
  [6,1,6; 6,7,6;11,11,-4;11,11,4;11,11,0],
  [8,5,-4;8,5,4],
  [6,1,6;6,7,6],
  [0,7,0],
  [3,7,0;6,1,6;6,7,6;3,7,4;3,7,-4],
  [8,5,8;8,5,4]
 };
 
end

ntonos = 2*length(K);

R = zeros(1,columns(M)-1);

% Constructs the consonant intervals

Ke = [];
for s = K;
  Ke = [[[0:ntonos-1];ones(1,ntonos)*s] Ke];
end

D = setdiff([0:ntonos-1],K);

De = [];
for s = D;
  De = [[[0:ntonos-1];ones(1,ntonos)*s] De];
end

for u=[0:ntonos-1]
 for v=[1:2:ntonos-1]
  if(((u!=0)||(v!=1))&&(sort(mod(K*v+u,ntonos))==D))
   plineal = v;
   pafin = u;
  end
 end
end

for w = [1:columns(M)-1]

% Finds the appropiate set of counterpoint
% symmetries.

 if(ismember(M(2,w),K))
  cual = find(K==M(2,w));
  esconsonancia = 1;
 else
  cual = find(K==mod(M(2,w)*plineal+pafin,ntonos));
  esconsonancia = 0;
 end
 
 for l = [1:rows(simetrias{cual})]
 
% Extracts symmetries as the matrix Q and the
% translation vector t.
 
  Q = [simetrias{cual}(l,2) 0;
  simetrias{cual}(l,3) simetrias{cual}(l,2)];
  t = [M(1,w); simetrias{cual}(l,1)];
  
% Calculation of the deformed intervals.

   consdeformadas = mod(Q*Ke+t*ones(1,length(Ke)),
   ntonos);
   disdeformadas = mod(plineal*Q*Ke+(t*plineal+
   pafin)*ones(1,length(Ke)),ntonos);   

% Adds one to the counter if the next
% interval is an admitted successor.
  if(esconsonancia&&ismember(M(2,w+1),K))
   if(ismember(M(:,w+1)',consdeformadas',"rows"))
    R(w) = R(w)+1;
   end
  elseif(esconsonancia&&ismember(M(2,w+1),D))
   if(ismember(M(:,w+1)',disdeformadas',"rows"))
    R(w) = R(w)+1;
   end
  elseif(not(esconsonancia)&&ismember(M(2,w+1),D))
   if(ismember(M(:,w+1)',disdeformadas',"rows"))
    R(w) = R(w)+1;
   end
  else
   if(ismember(M(:,w+1)',consdeformadas',"rows"))
    R(w) = R(w)+1;
   end  
  end
 end
end
end
\end{lstlisting}

The listed arrays contains the analyzed intervals extracted from Scriabin's sonata, one
for each part.

\begin{lstlisting}[language={Octave},basicstyle={\footnotesize}]
contrapunto_scriabin_13_31 = {
[4 4 4; 8 3 8],
[4 4 4;10 8 10],
[4 4 4; 6 10 6],
[4 4 4; 3 6 3],
[4 4; 10 3],
[4 4; 8 11],
[6 6; 10 3],
[6 6; 8 11],
[4 4 4; 2 1 0],
[8 8 8 8 8; 4 0 1 2 8],
[4 4 4; 2 6 2],
[8 8 8; 2 6 2],
[8 8 8; 3 6 3],
[8 8 8; 1 6 1],
[4 4 4; 8 7 6],
[10 10 10; 8 7 6],
[5 4 3 2 1; 11 11 11 11 11];
[1 4; 7 7],
[10 7; 8 5],
[4 1; 8 5]
};
contrapunto_scriabin_47_61 = {
[6 10 3;3 6 10],
[6 10 3; 11 3 6],
[3 6 3 11; 6 10 6 3],
[3 6 3 11; 11 3 10 8],
[3 11 10 11; 6 3 3 3],
[3 11 10 11; 10 8 6 8],
[3 11 3; 6 3 6],
[3 11 3; 10 8 10],
[1 4 1 9; 4 8 4 1],
[1 4 1 9; 9 1 8 6],
[8 9 8; 1 1 1],
[8 9 8; 4 6 4],
[7 10 7 3; 10 2 10 8],
[7 10 7 3; 3 5 3 0]
};	
\end{lstlisting}

\bibliographystyle{amsplain}
\bibliography{reporte18}

\end{document}